\documentclass[12pt]{article}
\usepackage{a4,latexsym,amsmath,amssymb,amsthm,enumerate,eucal,
  graphicx,subfigure}

\newtheorem{theorem}{Theorem}

\newtheorem{claim}{Claim}

\newcommand\Setx[1] {\left\{{#1}\right\}}
\newcommand\Set[2] {\left\{{#1}:\,{#2}\right\}}

\newcommand\size[1] {\left|{#1}\right|}

\newcommand{\PP}{\mathcal P}
\newcommand{\QQ}{\mathcal Q}
\newcommand{\TT}{\mathcal T}

\newcommand{\fig}[1]{\includegraphics[page=#1]{figs.pdf}}
\newcommand{\sfig}[2]{\begin{subfloat}\fig{#1}\caption{#2}\end{subfloat}}

\newcommand{\hf}{\hspace*{0mm}\hspace{\fill}\hspace*{0mm}}

\newbox\subfigbox
\makeatletter
\newenvironment{subfloat}
{\def\caption##1{\gdef\subcapsave{\relax##1}}
  \let\subcapsave=\@empty
  \let\sf@oldlabel=\label
  \def\label##1{\xdef\sublabsave{\noexpand\label{##1}}}
  \let\sublabsave\relax
  \setbox\subfigbox\hbox
  \bgroup}
{\egroup
  \let\label=\sf@oldlabel
  \subfigure[\subcapsave]{\box\subfigbox\sublabsave}}
\makeatother

\title{\textbf{A short proof\\of the tree-packing theorem}} %
\author{Tom\'{a}\v{s} Kaiser\footnote{Department of Mathematics and
    Institute for Theoretical Computer Science (ITI), University of
    West Bohemia, Univerzitn\'{\i}~8, 306~14~Plze\v{n}, Czech
    Republic. E-mail: \texttt{kaisert@kma.zcu.cz}. Supported by
    project 1M0545 and Research Plan MSM 4977751301 of the Czech
    Ministry of Education, and by project GA\v{C}R~201/09/0197 of the
    Czech Science Foundation.}  }

\date{}

\begin{document}
\maketitle

\begin{abstract}
  We give a short elementary proof of Tutte and Nash-Williams'
  characterization of graphs with $k$ edge-disjoint spanning trees.
\end{abstract}

We deal with graphs that may have parallel edges and loops; the vertex
and edge sets of a graph $H$ are denoted by $V(H)$ and $E(H)$,
respectively. Let $G$ be a graph. If $\PP$ is a partition of $V(G)$,
we let $G/\PP$ be the graph on the set $\PP$ with an edge joining
distinct vertices $X,Y\in\PP$ for every edge of $G$ with one end in
$X$ and another in $Y$. Tutte~\cite{Tut-problem} and
Nash-Williams~\cite{NW-edge} proved the following classical result:

\begin{theorem}\label{t:main}
  A graph $G$ contains $k$ pairwise edge-disjoint spanning trees if
  and only if for every partition $\PP$ of $V(G)$, the graph $G/\PP$
  has at least $k(\size\PP-1)$ edges.
\end{theorem}

Necessity of the condition in Theorem~\ref{t:main} is immediate. An
elegant proof of sufficiency is based on the matroid union theorem
(see, e.g., \cite[Corollary 51.1a]{Sch-combinatorial}) which yields
the more general matroid base packing theorem of
Edmonds~\cite{Edm-lehman}. A relatively short elementary proof of
sufficiency in Theorem~\ref{t:main}, due to W. Mader (personal
communication from R. Diestel), is given in~\cite[Theorem
2.4.1]{Die-graph}.

In this paper, we give another elementary proof that is also short and
perhaps somewhat more straightforward. The argument directly
translates to an efficient algorithm to find either $k$ disjoint
spanning trees, or a proof that none exist.

To give the reader an idea of the approach, let us briefly sketch the
proof of sufficiency, restricting to the case $k=2$. Let $T$ be a
spanning tree of $G$, and let $\overline T = G-E(T)$. We may assume
that $\overline T$ is disconnected as a spanning subgraph of $G$
(otherwise, we have two disjoint spanning trees). We seek a partition
$\PP$ of $V(G)$ such that each class of $\PP$ induces a connected
subgraph in both $T$ and $\overline T$. In order to find it, we start
with the trivial partition $\Setx{V(G)}$ and iteratively refine it (in
a suitable way) until we reach the desired partition $\PP$.

Let $E_\PP$ denote the set of edges of $G$ joining different classes
of $\PP$. The fact that $T[X]$ is connected for each $X\in\PP$ enables
us to count the edges of $T$ in $E_\PP$. Meanwhile, the density
condition yields a lower bound on $\size{E_\PP}$ and implies
$\size{E(\overline T) \cap E_\PP} \geq \size\PP-1$. Since $\overline
T$ is disconnected, and since $\overline T[X]$ is connected for all
$X\in\PP$, this forces a cycle in $\overline T$ intersecting at least
two classes of $\PP$. We can replace some edge of $T$ by an edge of
this cycle, so as to obtain a new spanning tree $T'$. When done
correctly, the exchange `improves' the spanning tree $T$ in a
well-defined way. Thus, if the initial spanning tree $T$ is chosen as
optimal, then the basic assumption that $\overline T$ is disconnected
must fail, which gives us the desired disjoint spanning trees.

A variant of this approach was used by Kaiser and
Vr\'{a}na~\cite{KV-hamilton} in connection with the conjecture of
Thomassen~\cite{Tho-reflections} that 4-connected line graphs are
hamiltonian. In that context, the method is applied to hypergraphs
instead of graphs and gives a connectivity condition under which a
hypergraph admits a `spanning hypertree' whose complement is, in a
way, close to being connected. A significant difference from the above
setup is that the situation in~\cite{KV-hamilton} is asymmetric
(unlike the packing of two spanning trees in a graph). It would be
interesting to identify more general conditions allowing for the
application of the method.

As noted by D. Kr\'{a}l' (personal communication), a matroid-theoretic
reformulation of the argument of the present paper yields a proof of
the matroid base packing theorem mentioned above.

Before we start with the detailed proof of Theorem~\ref{t:main}, we
introduce some terminology.  Let $k\geq 1$. A \emph{$k$-decomposition}
$\TT$ of a graph $G$ is a $k$-tuple $(T_1,\dots,T_k)$ of spanning
subgraphs of $G$ such that $\Set{E(T_i)}{1 \leq i \leq k}$ is a
partition of $E(G)$.

\begin{figure}
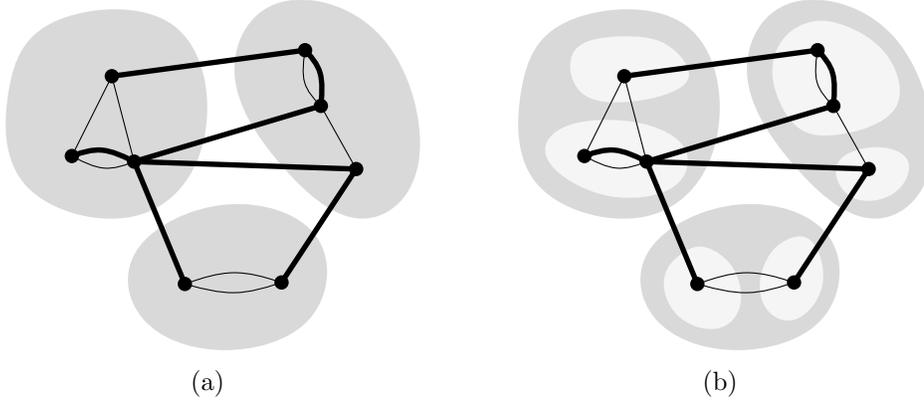

  \centering
  \hf
  \sfig1{}
  \hf
  \sfig2{}
  \hf
  \caption{The sequence of partitions associated with a 2-decomposition
    $\TT=(T_1,T_2)$ of $G$. The edges of $T_1$ are shown bold. (a) The
    partition $\PP_1$ (dark grey regions). (b) The partition $\PP_2$
    (light grey regions). Note that $\PP_2 = \PP_\infty$.}
  \label{fig:seq}
\end{figure}

We define the sequence $(\PP_0, \PP_1, \dots, \PP_\infty)$ of
partitions of $V(G)$ \emph{associated with $\TT$} as follows. (See the
illustration in Figure~\ref{fig:seq}.) First, $\PP_0 =
\Setx{V(G)}$. For $i\geq 0$, if there exists $c\in\Setx{1,\dots,k}$
such that the induced subgraph $T_c[X]$ is disconnected for some
$X\in\PP_i$, then let $c_i$ be the least such $c$, and let $\PP_{i+1}$
consist of the vertex sets of all components of $T_{c_i}[X]$, where
$X$ ranges over all the classes of $\PP_i$. Otherwise, the process
ends by setting $\PP_\infty = \PP_i$. In this case, we also set $c_j =
k+1$ and $\PP_j = \PP_i$ for all $j \geq i$.

The \emph{level} $\ell(e)$ of an edge $e\in E(G)$ (with respect to
$\TT$) is defined as the largest $i$ (possibly $\infty$) such that
both ends of $e$ are contained in one class of $\PP_i$. To keep the
notation simple, the symbols $\PP_i$ and $\ell(e)$ (as well as
$\PP_\infty$ and $c_i$) will relate to a $k$-decomposition $\TT$,
while $\PP'_i$ and $\ell'(e)$ will relate to a $k$-decomposition
$\TT'$. Thus, for instance, the level $\ell'(e)$ of an edge $e$ with
respect to $\TT'$ is defined using the partitions $\PP'_i$ associated
with $\TT'$.

When $\PP$ and $\QQ$ are partitions of $V(G)$, we say that $\PP$
\emph{refines} $\QQ$ (and write $\PP \leq \QQ$) if every class of
$\PP$ is a subset of a class of $\QQ$. When $\PP \leq \QQ$ and $\PP
\neq \QQ$, we write $\PP < \QQ$.

We define a strict partial order $\prec$ on $k$-decompositions of
$G$. Given two $k$-decompositions $\TT$ and $\TT'$, we set $\TT \prec
\TT'$ if there is some (finite) $j \geq 0$ such that both of the
following conditions hold:
\begin{enumerate}[(i)]
\item for $0 \leq i < j$, $\PP_i = \PP'_i$ and $c_i = c'_i$,
\item either $\PP_j < \PP'_j$, or $\PP_j = \PP'_j$ and $c_j < c'_j$.
\end{enumerate}

\begin{proof}[\textbf{Proof of Theorem~\ref{t:main}}]
  The necessity of the condition is clear. To prove the sufficiency,
  we proceed by induction on $k$. The claim is trivially true for
  $k=0$, so assume $k \geq 1$ and choose a $k$-decomposition
  $\TT=(T_1,\dots,T_k)$ of $G$ such that $T_1,\dots,T_{k-1}$ are trees
  and, subject to this condition, $\TT$ is maximal with respect to
  $\prec$.

  If $T_k$ is connected, then we are done. Otherwise, suppose that
  $T_k$ has at least two components (i.e., $\size{\PP_1}\geq 2$). We
  prove that there exists an edge of finite level (with respect to
  $\TT$) contained in a cycle of $T_k$. Let $\PP = \PP_\infty$.
  Recall that for $1\leq i < k$ and $X\in\PP$, the graph $T_i[X]$ is
  connected. Hence $T_i/\PP$ is a tree and has exactly $\size{\PP}-1$
  edges. By hypothesis, $G/\PP$ has at least $k(\size\PP-1)$ edges, so
  $T_k/\PP$ has at least $\size\PP-1$ edges. Since $T_k/\PP$ has
  $\size\PP$ vertices and is disconnected, it must contain a
  cycle. Thus $T_k$ contains a cycle, since $T_k[X]$ is connected for
  each $X\in\PP$. At least two edges of the cycle join different
  classes of $\PP$, and therefore their level is finite, as required.

  \begin{figure}
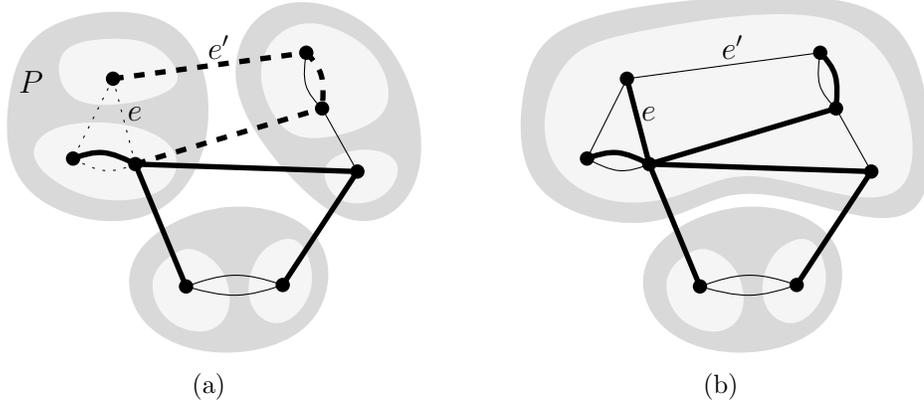

    \centering
    \hf
    \sfig3{}
    \hf
    \sfig4{}
    \hf
    \caption{The exchange step for the 2-decomposition $\TT$ of
      Figure~\ref{fig:seq}. (a) A cycle in $T_2$ containing $e$
      (dotted) and the cycle $C$ in $T_1+e$ (dashed). (b) The spanning
      tree $T'_1$ (bold) obtained from $T_1$ by exchanging $e$ for the
      edge $e'$ of $C$. The partitions $\PP'_1$ and $\PP'_2$
      associated with the resulting 2-decomposition $\TT'$ are shown
      in dark grey and light grey, respectively. Note that $\PP'_2$ is
      equal to $\PP'_\infty$ and that $\TT \prec \TT'$.}
    \label{fig:change}
  \end{figure}
  
  Let $e\in E(T_k)$ be an edge of minimum level that is contained in a
  cycle of $T_k$, and set $m=\ell(e)$. (See Figure~\ref{fig:change}
  for an illustration with $m=1$.) Let $P$ be the class of $\PP_m$
  containing both ends of $e$. Since $e$ joins different components of
  $T_{c_m}[P]$, we have $c_m\neq k$, and the unique cycle $C$ in
  $T_{c_m}+e$ contains an edge with only one end in $P$. Thus, for an
  edge $e'$ of $C$ of lowest possible level we have $\ell(e') <
  m$. Let $Q$ be the class of $\PP_{\ell(e')}$ containing both ends of
  $e'$. Observe that $V(C) \subseteq Q$. We will exchange $e$ for $e'$
  in the members of the $k$-decomposition to eventually obtain the
  desired contradiction.

  Let $\TT'$ be the $k$-decomposition obtained from $\TT$ by replacing
  $T_{c_m}$ with $T_{c_m} + e - e'$ and $T_k$ with $T_k - e + e'$. The
  $i$-th element of $\TT'$, where $1\leq i \leq k$, is denoted by
  $T'_i$. To relate the sequences of partitions associated with $\TT$
  and $\TT'$, we prove the following two claims.

  \begin{claim}\label{cl:connected}
    If $T_c[X]$ is connected, for some $X\subseteq V(G)$ and $1\leq c
    \leq k$, then $T'_c[X]$ is connected unless one of the following
    holds:
    \begin{enumerate}[\quad(a)]
    \item $c = c_m$, and $X$ contains both ends of $e'$, and
      $Q\not\subseteq X$, or
    \item $c = k$, and $X$ contains both ends of $e$, and
      $P\not\subseteq X$.
    \end{enumerate}
  \end{claim}
  
  To prove the claim, suppose that $T'_c[X]$ is disconnected. We have
  $c\in\Setx{c_m,k}$, since otherwise $T_c = T'_c$. Consider $c =
  c_m$. Since $E(T_{c_m}) - E(T'_{c_m}) = \Setx {e'}$, both ends of
  $e'$ lie in $X$. Furthermore, $Q\not\subseteq X$, since otherwise
  $T'_{c_m}[X]$ would contain the path $C-e'$ joining the ends of
  $e'$, which would make $T'_{c_m}[X]$ connected. A similar argument
  for the case $c=k$ completes the proof of Claim~\ref{cl:connected}.

  \begin{claim}\label{cl:same}
    For all $i \leq m$, it holds that $c'_i = c_i$ and $\PP'_i =
    \PP_i$.
  \end{claim}
  We proceed by induction on $i$. The case $i=0$ follows from $\PP_0 =
  \PP'_0 = \Setx{V(G)}$ and $c_0 = c'_0 = k$. Let us thus assume that
  the assertion holds for some $i$, $0\leq i < m$, and prove it for
  $i+1$.

  We first prove that $\PP_{i+1} = \PP'_{i+1}$. Let $S$ be an
  arbitrary class of $\PP_{i+1}$; we assert that $T'_{c'_i}[S]$ is
  connected. Since $T_{c_i}[S]$ is connected and since $c'_i = c_i$ by
  the inductive hypothesis, we can use Claim~\ref{cl:connected} (with
  $X=S$ and $c=c_i$). Condition (a) in the claim cannot hold, because
  every class of $\PP_{i+1}$ containing both ends of $e'$ contains $Q$
  as a subset. For a similar reason, condition (b)
  fails. Consequently, $T'_{c_i}[S]$ is connected, and hence $S$ is a
  subset of some class of $\PP'_{i+1}$. Since $S$ was arbitrary, it
  follows that $\PP_{i+1} \leq \PP'_{i+1}$. Now by the choice of $\TT$
  (and the inductive assumption), we cannot have $\PP_{i+1} <
  \PP'_{i+1}$. We conclude that $\PP_{i+1} = \PP'_{i+1}$.

  Next, we prove that $c'_{i+1} = c_{i+1}$. Let $R\in\PP'_{i+1}$ and
  $c < c_{i+1}$. By the above, $R\in\PP_{i+1}$. The definition of
  $c_{i+1}$ implies that $T_c[R]$ is connected. Using
  Claim~\ref{cl:connected} as above, we find that $T'_c[R]$ is also
  connected. Consequently, $c'_{i+1} \geq c_{i+1}$, and by the
  maximality of $\TT$ once again, we must have $c'_{i+1} =
  c_{i+1}$. The proof of Claim~\ref{cl:same} is complete.

  It is now easy to finish the proof of Theorem~\ref{t:main}. Since
  $\PP'_m = \PP_m$ and $c'_m = c_m$, the classes of $\PP'_{m+1}$ are
  the vertex sets of components of $T'_{c_m}[U]$, where
  $U\in\PP_m$. Observe that for $U\in\PP_m-\Setx P$, we have
  $T'_{c_m}[U] = T_{c_m}[U]$, and so the components of $T'_{c_m}[U]$
  coincide with those of $T_{c_m}[U]$. The graph $T'_{c_m}[P]$ is
  obtained from $T_{c_m}[P]$ by adding the edge $e$ that connects two
  components of $T_{c_m}[P]$. It follows that $\PP_{m+1} <
  \PP'_{m+1}$, contradicting the choice of $\TT$.
\end{proof}

\section*{Acknowledgment}
\label{sec:acknowledgment}

I am indebted to Douglas West and two anonymous referees who suggested
a number of improvements to the paper.

\end{document}